\documentclass[11pt]{article}
\usepackage{amsmath, amsfonts, amssymb, amsthm}
\usepackage{graphicx, mathrsfs, dsfont}
\usepackage{enumitem}
\usepackage{comment}
\usepackage[numbers,sort&compress]{natbib}
\usepackage{hyperref}
\let\svthefootnote\thefootnote
\newcommand\freefootnote[1]{%
  \let\thefootnote\relax%
  \footnotetext{#1}%
  \let\thefootnote\svthefootnote%
}
\usepackage{tikz,pgf}
\usepackage{pgfplots}
\usetikzlibrary{shapes.geometric}
\usetikzlibrary{arrows}
\usetikzlibrary{automata,patterns,calc}
\usepackage{tikz-3dplot}
\pgfplotsset{compat=1.18}
\usepackage[T1]{fontenc}

\usepackage{float}
\providecommand{\U}[1]{\protect\rule{.1in}{.1in}}

\newtheorem{theorem}{Theorem}[section]
\newtheorem{corollary}[theorem]{Corollary}
\newtheorem{lemma}[theorem]{Lemma}
\newtheorem{proposition}[theorem]{Proposition}
\newtheorem{assumption}{Assumption}

\theoremstyle{definition}
\newtheorem{remark}[theorem]{Remark}

\newtheorem{example}[theorem]{Example}

\providecommand{\gph}{\mathop{\rm gph}\nolimits}

\providecommand{\into}{\mathop{\rm int}\nolimits}

\providecommand{\argmax}{\mathop{\rm argmax}\nolimits}

\providecommand{\tto}{\mathop{\rightrightarrows}\nolimits}

\newcommand{\R}{\mathbb{R}}

\newcommand{\N}{\mathbb{N}}

\title{Applicability of the MPCC approach to Zeng reformulations of nonlinear pessimistic bilevel programming}
\author{ Joaqu\'in Parra, David Salas}
\date{}

\begin{document}

\maketitle

\begin{abstract}
The vanilla method in bilevel programming is to reformulate it as a single-level mathematical programs with complementarity constraints (MPCC), via the KKT conditions of the lower-level. This method seemed to work only for optimistic formulations. However, in \cite{Zeng:2020}, B. Zeng presented a general method to reduce pessimistic bilevel problems to optimistic ones. In this note, we study how to apply the MPCC approach to the pessimistic counterpart through Zeng's reformulation \cite{Zeng:2020}. While we show that it cannot be applied directly since constraint qualification can be destroyed by the reformulation, it still can be used to solve natural relaxations, in the spirit of \cite{Wiesemann:2013,Loridan-Morgan:1989,LoridanMorgan:1996}.


\end{abstract}

\section{Introduction}

Bilevel programming models a leader who anticipates the optimal response of a follower. When multiple optimal responses are possible, optimistic bilevel programming considers the most favorable response for the leader, whereas pessimistic bilevel programming considers the least favorable one. Its min--max--min structure introduces an additional worst-case optimization step, which makes the pessimistic formulation appear more difficult even in the linear setting. 

In a very unexpected result, B. Zeng \cite{Zeng:2020} showed a simple way to reduce a pessimistic formulation to an optimistic one. Although the computational development in \cite{Zeng:2020} focuses on the linear setting, the reformulation also applies to nonlinear problems under standard compactness and continuity assumptions.

Zeng's reformulation \cite{Zeng:2020} has contributed to our understanding of bilevel programming, for example in the study of coupling constraints \cite{Henke-et-al:2025:Coupling,Henke-et-al:2026:PessimisticCoupling}, and in identifying the gap between bilevel programming and single-leader multi-follower games \cite{AusselSalasSchmidt:2026}. However, its scope for algorithmic applications remains to be further explored.

In optimistic bilevel programming, a standard approach is to replace the lower-level optimality requirement by the Karush--Kuhn--Tucker (KKT) conditions, under suitable convexity and constraint-qualification assumptions \cite{Dempe:2002}. The resulting single-level problem is a mathematical program with complementarity constraints (MPCC), for which there is an extensive theoretical and numerical literature; see, e.g., \cite{HoheiselKanzowSchwartz:2013,KanzowSchwartz:2013} and the references therein. In the pessimistic setting, this direct substitution does not eliminate the worst-case maximization, which has motivated alternative approaches. 

Existing work in this vein includes optimality conditions \cite{DempeEtAl:2014} and relaxation methods \cite{Zemkoho1:2025,Zemkoho2:2026}.
%, and proximal and gradient methods \cite{GuanEtAl:2025,CaoZengZhang:2025}. 
Specialized computational approaches have also been developed for problems with linear or discrete structure \cite{GoyalZhangHe:2023,CalveteEtAl:2025,BucareyEtAl:2026}. More recently, Stein and Zemkoho \cite{Stein-Zemkoho:2026} proposed a single-level reformulation that does not use complementarity constraints or optimal value functions, offering another approach to the pessimistic problem.

In this note, we study the applicability of the MPCC approach to Zeng reformulations in the nonlinear setting. We show that constraint qualification may fail for the exact reformulation, but is preserved by a positive relaxation of its additional sublevel constraint, in the spirit of \cite{Loridan-Morgan:1989,LoridanMorgan:1996,Wiesemann:2013}. Since in general the relaxed problems may not attain their infima, we prove a convergence result for $\delta_n$-minimizers as the relaxation and optimality tolerances tend to zero. This establishes a basis for applying the MPCC approach to pessimistic bilevel programming via Zeng's reformulation. We also discuss an extension to robust pessimistic constraints through penalization.

\subsection{Preliminaries}

We write $x\in\R^p$ for the leader's decision variable and $y\in\R^q$ for the follower's decision variable. For $n\in\N$, we use $[n]=\{1,\ldots,n\}$, and for $a,b\in\R^n$, the product $a\odot b=(a_ib_i)_{i\in[n]}$ is understood componentwise. Vector inequalities are also componentwise, and $\R_+=[0,+\infty)$. For a set-valued map $M:D\tto\R^q$, its graph is denoted by $\gph M=\{(u,v)\in D\times\R^q\mid v\in M(u)\}$. All continuity and semicontinuity properties are relative to the stated domains.

An exact solution of a minimization problem is called a \emph{global minimizer}. For $\delta\geq 0$, a \emph{$\delta$-minimizer} of $\min_{z\in\Omega}f(z)$ is a feasible point $z\in\Omega$ such that $f(z)\leq\inf_{w\in\Omega}f(w)+\delta$. Thus, $0$-minimizers are global minimizers. The verb ``solve'' always refers to exact global optimality, for both minimization and maximization problems. In a bilevel problem, the tolerance $\delta$ applies only to the upper-level objective: feasibility still requires the follower to solve the lower-level problem exactly.

Our MPCC reformulations use differentiability with respect to the follower's variable. The assumptions needed for these reformulations and those used for the approximation results are stated separately below.

\paragraph{Constraint qualification for convex optimization problems:} Let us consider an abstract optimization problem
\begin{equation}\label{eq:Problemonelevel}
    \min_{z\in \Omega} f(z),
\end{equation}
where the feasible set $\Omega\subset \R^{n}$ is nonempty, convex and closed, and is defined by the constraints
\begin{equation}\label{eq:General-constraintset}
    \Omega = \left\{z\in\R^n\left| \begin{array}{l}
        h_i(z) = 0,\quad \forall i\in I,\\
        g_j(z) \leq 0,\quad \forall j\in J
    \end{array}\right.\right\},
\end{equation}
where $I$ and $J$ are finite index sets and the functions $h_i,g_j:\R^{n}\to\R$ are differentiable. For $z\in\Omega$, we define the \textit{normal cone} of $\Omega$ at $z$ (in the sense of convex analysis) by
\begin{equation}\label{eq:normalConeDef}
N_{\Omega}(z)= \{ \nu\in\R^n~:~\langle \nu,y-z\rangle \leq 0,\,\forall y\in\Omega \}.
\end{equation}
Similarly, we define the \emph{Farkas cone} associated with the representation \eqref{eq:General-constraintset} by
\begin{equation}\label{eq:FarkasCone}
    \mathcal{F}_{\Omega}(z) = \left\{ \sum_{i\in I} \lambda_i\nabla h_i(z) + \sum_{j\in J} \mu_j\nabla g_j(z)\, \Big|\ \mu\odot g(z) = 0,\; \mu\geq 0\right\}.
\end{equation}
Here $\lambda\in\R^I$, $\mu\in\R^J$, and $g=(g_j)_{j\in J}$. We will say that the set $\Omega$, with its representation \eqref{eq:General-constraintset}, is \emph{constraint qualified} at $z\in \Omega$ if
\begin{equation}\label{eq:CQ}
    N_{\Omega}(z) = \mathcal{F}_{\Omega}(z).
\end{equation}
We say that the representation is constraint qualified if this equality holds at every point of $\Omega$.
Some works refer to condition \eqref{eq:CQ} as Guignard's constraint qualification, as presented in \cite{guignard1969generalized} in arbitrary Banach spaces. For a differentiable convex objective, condition \eqref{eq:CQ} allows us to characterize global minimizers by their KKT conditions. There are many sufficient conditions for \eqref{eq:CQ}. In particular, if all the functions $h_i$ and $g_j$ are affine, then \eqref{eq:CQ} holds (see, e.g., \cite{solodov2010constraint}).

\paragraph{MPCC reformulations of optimistic problems:}

Consider the optimistic bilevel programming problem
\begin{equation}\label{eq:Optimistic}
	\begin{array}{cl}
		\displaystyle\min_{x,y} & F(x,y)\\
		\text{s.t.} &\begin{cases}
			x\in X,\\
			y\text{ solves }\left\{\begin{array}{cl}
				\displaystyle\min_y& f(x,y)  \\
				\text{s.t.} & \begin{cases}
					h_i(x,y) = 0,\quad&\forall i\in I,\\
                    g_j(x,y) \leq 0,\quad&\forall j\in J.
				\end{cases}
			\end{array}\right.
		\end{cases}
	\end{array}
\end{equation}
Suppose that, for every $x\in X$, the functions $f(x,\cdot)$, $h_i(x,\cdot)$ and $g_j(x,\cdot)$ are differentiable on $\R^q$. The following two hypotheses are sufficient to obtain an MPCC reformulation of Problem~\eqref{eq:Optimistic}:
\begin{enumerate}
    \item \textbf{Convexity of the problem:} For every $x\in X$, the functions $f(x,\cdot)$ and $g_j(x,\cdot)$ are convex, and the functions $h_i(x,\cdot)$ are affine.
    \item \textbf{Constraint qualification of the problem:} For every $x\in X$, the feasible set of the lower-level problem, with the representation in \eqref{eq:Optimistic}, is constraint qualified in the sense of \eqref{eq:CQ}.
\end{enumerate}
Under these hypotheses, Problem~\eqref{eq:Optimistic} is equivalent to
\begin{equation}\label{eq:MPCC}
	\begin{array}{cl}
		\displaystyle\min_{x,y,\lambda,\mu} & F(x,y)\\
		\text{s.t.} &\begin{cases}
			x\in X,\\
			h_i(x,y) = 0,\quad&\forall i\in I,\\
            g_j(x,y) \leq 0,\quad&\forall j\in J,\\
            \displaystyle\nabla_yf(x,y) + \sum_{i\in I} \lambda_i\nabla_y h_i(x,y) + \sum_{j\in J} \mu_j\nabla_y g_j(x,y) = 0,\\
			\mu\odot g(x,y) = 0,\\
            \mu\geq 0.
		\end{cases}
	\end{array}
\end{equation}
More precisely, the equivalence is understood as
\[
(x,y)\text{ solves~\eqref{eq:Optimistic}}\iff\exists(\lambda,\mu)\text{ such that }(x,y,\lambda,\mu)\text{ solves~\eqref{eq:MPCC}}.
\]
%Problem~\eqref{eq:MPCC} is the MPCC reformulation of Problem~\eqref{eq:Optimistic}. The equivalence above concerns global minimizers; it does not assert a correspondence between local minimizers. For a pessimistic problem, replacing the original follower's problem by its KKT conditions does not by itself eliminate the maximization over its solutions. We study the use of this construction after the Zeng reformulation.

\section{Pessimistic formulations and assumptions}\label{sec:setting}

We consider a nonempty set $X\subset\R^p$ and, for every $x\in X$, the parametric optimization problem
\begin{equation}\label{eq:lower-level}
    \begin{array}{cl}
         \displaystyle\min_{y\in \R^q}& f(x,y)  \\
         \text{s.t.} & \begin{cases}
             h_i(x,y) = 0,\quad \forall i\in I,\\
             g_j(x,y) \leq 0,\quad \forall j\in J,
         \end{cases}
    \end{array}
\end{equation}
where $I$ and $J$ are finite index sets and the functions $f,h_i,g_j$ are defined on $X\times\R^q$. In what follows, we write
\begin{equation}\label{eq:Feasible-lower-level}
    K(x) = \left\{y \in \R^q \left| \begin{array}{ll}
        h_i(x,y) = 0,& \forall i\in I,\\
        g_j(x,y) \leq 0,& \forall j\in J
    \end{array}  \right.\right\}
\end{equation}
to denote the feasible set of Problem~\eqref{eq:lower-level}. We assume that $K(x)$ is nonempty for every $x\in X$ and that there exists a compact set $Y\subset\R^q$ such that $K(x)\subset Y$ for all $x\in X$. We also assume that the restrictions of $f,h_i,g_j$ to $X\times Y$ are continuous. Thus, the solution set of~\eqref{eq:lower-level} is nonempty and compact for every $x\in X$.

Denoting by $S:X\tto Y$ the solution mapping of Problem~\eqref{eq:lower-level}, we consider the pessimistic bilevel problem
\begin{equation}\label{eq:Pessimistic}
    \min_{x\in X}\max_{y\in S(x)} F(x,y),\tag{PBP}
\end{equation}
where $F:X\times\R^q\to\R$ is also continuous on $X\times Y$. We adopt a \emph{leader-only} notion of solution for~\eqref{eq:Pessimistic}: a global minimizer is a point $x^*\in X$ that globally minimizes $\varphi:x\mapsto\max_{y\in S(x)}F(x,y)$. We summarize the standing compactness and continuity assumptions below.
\begin{assumption}\label{as:standard} We consider the following standard assumptions:
\begin{enumerate}
    \item[(i)] $X\subset\R^p$ and $Y\subset\R^q$ are nonempty compact sets.
    \item[(ii)] The functions $F,f,h_i,g_j:X\times\R^q\to\R$ are continuous on $X\times Y$.
    \item[(iii)] For every $x\in X$, the set $K(x)$ defined in~\eqref{eq:Feasible-lower-level} is nonempty and contained in $Y$.
\end{enumerate}
\end{assumption}

Assumption~\ref{as:standard} entails that $S$ takes nonempty compact values, so Problem~\eqref{eq:Pessimistic} is well-defined. To obtain the MPCC reformulations studied in this work, we impose additional differentiability, convexity and constraint qualification assumptions. In particular, we require $F(x,\cdot)$ to be concave in order to use convex optimality conditions for the inner maximization.
\begin{assumption}\label{as:convexity} We consider the following convexity-CQ assumptions on Problem~\eqref{eq:Pessimistic}:
	\begin{enumerate}
		\item[(i)] For every $x\in X$, the function $F(x,\cdot):\R^q\to\R$ is concave and differentiable.
		\item[(ii)] For every $x\in X$, the functions $f(x,\cdot),g_j(x,\cdot):\R^q\to\R$ are convex and differentiable, and the functions $h_i(x,\cdot):\R^q\to\R$ are affine.
		\item[(iii)] For every $x\in X$, the set $K(x)$, with its representation in~\eqref{eq:Feasible-lower-level}, is constraint qualified.
	\end{enumerate}
\end{assumption}

In general, Problem~\eqref{eq:Pessimistic} may fail to have global minimizers even under Assumptions~\ref{as:standard} and~\ref{as:convexity}. A standard sufficient condition for existence is lower semicontinuity of the function
\begin{equation}\label{eq:Pessimistic-function}
\begin{aligned}
\varphi:X&\to\R\\
x&\mapsto \varphi(x):= \max_{y\in S(x)}F(x,y).
\end{aligned}
\end{equation}
In this case, compactness of $X$ and the Weierstrass theorem ensure the existence of global minimizers and imply that all accumulation points of minimizing sequences are global minimizers. Much of the literature on existence either provides sufficient conditions for lower semicontinuity of $\varphi$ (usually through lower semicontinuity of the set-valued map $S$; see, e.g., \cite{Wiesemann:2013,LamparielloEtAl:2019,LiuEtAl:2018} and the references therein), or studies particular cases where this property can be verified directly.

\begin{assumption}\label{as:main} We assume that the function $\varphi$ defined in~\eqref{eq:Pessimistic-function} is lower semicontinuous.
\end{assumption}

\section{Zeng reformulation}\label{sec:zeng}

The Zeng reformulation of Problem~\eqref{eq:Pessimistic}, introduced in \cite{Zeng:2020}, is the following optimistic bilevel problem:
\begin{equation}\label{eq:Zeng-0}
\tag{$\mathcal{Z}$}
	\left.\begin{array}{cl}
		\displaystyle\min_{x,\hat{y},y} & F(x,y)\\
		\text{s.t.} &\begin{cases}
			x\in X,\ \hat{y}\in K(x),\\
			y\text{ solves }\left\{\begin{array}{cl}
				\displaystyle\max_y& F(x,y)  \\
				\text{s.t.} & \begin{cases}
					y\in K(x),\\
					f(x,y)\leq f(x,\hat{y}).
				\end{cases}
			\end{array}\right.
		\end{cases}
	\end{array}\right.
\end{equation}
Here $(x,\hat{y})$ are the leader's variables and $y$ is the follower's variable. Under Assumptions~\ref{as:standard} and~\ref{as:convexity}, the lower-level problem is convex.

Zeng \cite{Zeng:2020} calls~\eqref{eq:Zeng-0} a ``relaxation'' when solutions are represented by leader--follower pairs $(x,y)$. With the leader-only notion of solution adopted here, the two problems are equivalent in terms of global minimizers. We provide a short proof for completeness.

\begin{proposition}\label{prop:zeng-equivalence} Under Assumption~\ref{as:standard}, Problems~\eqref{eq:Pessimistic} and~\eqref{eq:Zeng-0} are equivalent in the sense that
	\[
	x\text{ solves~\eqref{eq:Pessimistic}}\iff \exists \hat{y},y\in Y\text{ such that } (x,\hat{y},y)\text{ solves~\eqref{eq:Zeng-0}}.
	\]
\end{proposition}
\begin{proof}
	Let $x^*$ solve~\eqref{eq:Pessimistic} and choose $y^*\in S(x^*)$ such that $F(x^*,y^*)=\varphi(x^*)$. Then $(x^*,y^*,y^*)$ is feasible for~\eqref{eq:Zeng-0}. For any other feasible point $(x,\hat{z},z)$, the inequality $f(x,\hat{z})\geq\min_{v\in K(x)}f(x,v)$ implies that every point of $S(x)$ is feasible for its lower-level problem. Hence,
	\[
	F(x,z)\geq\varphi(x)\geq\varphi(x^*)=F(x^*,y^*).
	\]
	Thus, $(x^*,y^*,y^*)$ solves~\eqref{eq:Zeng-0}.

	Conversely, let $(x^*,\hat{y}^*,y^*)$ solve~\eqref{eq:Zeng-0}. Feasibility gives $F(x^*,y^*)\geq\varphi(x^*)$. Choose $z^*\in S(x^*)$ attaining $\varphi(x^*)$. Comparing with the feasible point $(x^*,z^*,z^*)$ yields the reverse inequality, so $F(x^*,y^*)=\varphi(x^*)$. Now, fix $x\in X$ and choose $y\in S(x)$ with $F(x,y)=\varphi(x)$. Since $(x,y,y)$ is feasible for~\eqref{eq:Zeng-0},
	\[
	\varphi(x^*)=F(x^*,y^*)\leq F(x,y)=\varphi(x).
	\]
	Therefore, $x^*$ solves~\eqref{eq:Pessimistic}.
\end{proof}

\section{Preservation of constraint qualifications}\label{sec:Preservation}

To apply an MPCC reformulation to Problem~\eqref{eq:Zeng-0}, we need its lower-level problem to be convex and constraint qualified. Assumption~\ref{as:convexity} ensures convexity but we also need that the set-valued map $\widehat{K}:\gph K \tto Y$ given by
    \begin{equation}\label{eq:Feasible-Z0}
    \widehat{K}(x,\hat{y}) = \{y\in \R^q\mid y\in K(x),\ f(x,y)\leq f(x,\hat{y})\},
    \end{equation}
    to be constraint qualified in the sense of~\eqref{eq:CQ}, for every $x\in X$ and $\hat{y}\in K(x)$. This latter condition may fail even when the original representation of $K(x)$ is constraint qualified. The following example, adapted from \cite[Example~2.2]{DempeDutta2012}, illustrates this fact.

\begin{example}\label{ex:loss-CQ} Consider Problem~\eqref{eq:Pessimistic} with the following data:
\begin{itemize}
\item $X=\{0\}$ and $Y=[-1,1]^2\subset\R^2$.
\item $F(x,y)=-y_1$.
\item $f(x,y)=y_1^2+y_2$.
\item $K(x):=\left\{y\in\R^2\left|\begin{array}{c}-1\leq y_1,y_2\leq 1,\\
y_1^2-y_2\leq 0\end{array}\right.\right\}$.
\end{itemize}
The data satisfy Assumption~\ref{as:standard} and, since $X$ is a singleton, Assumption~\ref{as:main}. Moreover, the lower-level problem is convex and the representation of $K(0)$ satisfies Slater's condition, so Assumption~\ref{as:convexity} also holds.

The Zeng reformulation is
\begin{equation}\label{eq:ex-Z0}
\begin{array}{cl}
    \displaystyle\min_{(x,\hat{y},y)} & -y_1\\
    \text{s.t.} & x=0,\ \hat{y}\in K(0),\\
    &y\text{ solves }
    \left\{
    \begin{array}{ll}
        \min_y & y_1\\
        \text{s.t.} & \begin{cases}
        -1\leq y_1,y_2\leq 1,\\
        y_1^2-y_2\leq 0,\\
        y_1^2+y_2\leq\hat{y}_1^2+\hat{y}_2.
        \end{cases}
    \end{array}
    \right.
\end{array}
\end{equation}
Its unique global minimizer is $(x^*,\hat{y}^*,y^*)=(0,(0,0),(0,0))$. However, the representation of the lower-level feasible set at $(0,\hat{y}^*)$ is not constraint qualified. Indeed, $\widehat{K}(0,\hat{y}^*)$ is a singleton, so
\[
N_{\widehat{K}(0,\hat{y}^*)}(y^*)=\R^2.
\]
The gradients of the two active constraints at $y^*=(0,0)$ are $(0,-1)$ and $(0,1)$; the box constraints are inactive. Thus,
\[
\mathcal{F}_{\widehat{K}(0,\hat{y}^*)}(y^*)=\R e_2.
\]
In particular, these gradients cannot balance the horizontal gradient of the lower-level objective, so no choice of multipliers makes this solution feasible for the corresponding MPCC reformulation.\hfill$\Diamond$
\end{example}

The difficulty is that the additional constraint $f(x,y)\leq f(x,\hat{y})$ restricts $K(x)$ to $S(x)$ whenever $\hat{y}\in S(x)$. This may reduce its dimension and destroy constraint qualification. Following the discussion in \cite{Zeng:2020}, and inspired by classical regularization approaches (see, e.g., \cite{Wiesemann:2013,LoridanMorgan:1996,Loridan-Morgan:1989} or the survey \cite{CarusoLignolaMorgan:2020}), we consider, for $\varepsilon>0$, the relaxed problem
\begin{equation}\label{eq:Zeng-eps}
\tag{$\mathcal{Z}_{\varepsilon}$}
\left.\begin{array}{cl}
\displaystyle\min_{x,\hat{y},y} & F(x,y)\\
\text{s.t.} &\begin{cases}
    x\in X,\ \hat{y}\in K(x),\\
    y\text{ solves }\left\{\begin{array}{cl}
         \displaystyle\max_y& F(x,y)\\
         \text{s.t.} & \begin{cases}
             y\in K(x),\\
             f(x,y)\leq f(x,\hat{y})+\varepsilon.
         \end{cases}
    \end{array}\right.
    \end{cases}
\end{array}\right.
\end{equation}

\begin{proposition}\label{prop:preservation-CQ} For $x\in X$, $\hat{y}\in K(x)$ and $\varepsilon\geq 0$, consider the set
\begin{equation}\label{eq:widehatK}
    \widehat{K}(x,\hat{y},\varepsilon):=\left\{y\in\R^q\left|\begin{array}{l}
         h_i(x,y)=0,\quad\forall i\in I,\\
         g_j(x,y)\leq 0,\quad\forall j\in J,\\
         f(x,y)\leq f(x,\hat{y})+\varepsilon
    \end{array}\right.\right\}.
\end{equation}
Suppose that $K(x)$ is closed and convex, its representation in~\eqref{eq:Feasible-lower-level} is constraint qualified, and the functions $h_i(x,\cdot)$ and $g_j(x,\cdot)$ are differentiable on $\R^q$. If $f(x,\cdot):\R^q\to\R$ is convex and differentiable and $\varepsilon>0$, then the representation~\eqref{eq:widehatK} is also constraint qualified.
\end{proposition}
\begin{proof}
Let $C=\{y\in\R^q\mid f(x,y)\leq f(x,\hat{y})+\varepsilon\}$. Since $\varepsilon>0$ and $f(x,\cdot)$ is continuous, $\hat{y}\in\into C$. Thus, $K(x)\cap\into C\neq\emptyset$. The convex normal-cone intersection rule (see, e.g., \cite[Theorem~6.42]{Rockafellar-Wets:1998}) therefore gives, for every $y\in K(x)\cap C$,
\[
N_{K(x)\cap C}(y)=N_{K(x)}(y)+N_C(y)
=\mathcal{F}_{K(x)}(y)+\mathcal{F}_C(y).
\]
The last equality follows from constraint qualification of the representation of $K(x)$ and Slater's condition for the sublevel constraint defining $C$. Since $K(x)\cap C=\widehat{K}(x,\hat{y},\varepsilon)$, this identity means that $\nu\in N_{\widehat{K}(x,\hat{y},\varepsilon)}(y)$ if and only if there exist $\lambda\in\R^I$, $\mu\in\R^J$ and $\eta\in\R$ such that
\begin{gather*}
\nu=\sum_{i\in I}\lambda_i\nabla_yh_i(x,y)
+\sum_{j\in J}\mu_j\nabla_yg_j(x,y)+\eta\nabla_yf(x,y),\\
\mu\geq 0,\quad\eta\geq 0,\quad\mu\odot g(x,y)=0,\quad
\eta\bigl(f(x,y)-f(x,\hat{y})-\varepsilon\bigr)=0.
\end{gather*}
This is precisely the Farkas representation $\mathcal{F}_{\widehat{K}(x,\hat{y},\varepsilon)}(y)$ for~\eqref{eq:widehatK} at $y$, which proves the result.
\end{proof}

For $\varepsilon=0$, the notation in~\eqref{eq:widehatK} agrees with~\eqref{eq:Feasible-Z0}, that is, $\widehat{K}(x,\hat{y},0)=\widehat{K}(x,\hat{y})$. We obtain the following direct corollary.
\begin{corollary}\label{cor:MPCC-relaxation} Consider Assumptions~\ref{as:standard} and~\ref{as:convexity}.
Then, for every $\varepsilon>0$, Problem~\eqref{eq:Zeng-eps} is equivalent to its MPCC reformulation.
\end{corollary}

Equivalence to an MPCC does not guarantee that a relaxed problem has global minimizers, even when the exact problem does. To discuss this issue, under Assumption~\ref{as:standard}, let us define the inner optimal value function $\widehat{\varphi}:\gph K\times\R_+\to \R$ given by
\[
\widehat{\varphi}(x,\hat{y},\varepsilon)
:=\max_{y\in\widehat{K}(x,\hat{y},\varepsilon)}F(x,y).
\]
The following example illustrates the possible failure of attainment.

\begin{example}\label{ex:nonattainment}
Assumptions~\ref{as:standard},~\ref{as:convexity} and~\ref{as:main} do not guarantee
that Problem~\eqref{eq:Zeng-eps} has a global minimizer, even when
the lower-level objective is strongly convex and the solution
mapping $S$ is constant. Consider
\[
X=Y=[0,1],\qquad
f(x,y)=y^2,\qquad
F(x,y)=x+y,
\]
and
\[
K(x)=\{y\in\R\mid xy=0,\ 0\leq y\leq 1\}.
\]
All the functions are polynomial, and $K(x)$ is
nonempty and contained in $Y$ for every $x\in X$.
Thus, Assumption~\ref{as:standard} holds.
Moreover, for every fixed $x\in X$, the function $f(x,\cdot)$
is strongly convex, $F(x,\cdot)$ is affine, and all the
constraints defining $K(x)$ are affine with respect to $y$.
In particular, this representation of $K(x)$ is constraint
qualified, so Assumption~\ref{as:convexity} also holds.

We have
\[
K(0)=[0,1],
\qquad
K(x)=\{0\}\quad\text{for }x>0.
\]
Since $f(x,y)=y^2$, it follows that
\[
S(x)=\{0\}
\quad\text{and}\quad
\varphi(x)=F(x,0)=x
\qquad\text{for every }x\in X.
\]
Consequently, $\varphi$ is continuous, so
Assumption~\ref{as:main} holds, and the original pessimistic
problem has the unique global minimizer $x=0$.

Now fix $0<\varepsilon<1$ and consider Problem~\eqref{eq:Zeng-eps}.
If $x>0$, then $\hat{y}\in K(x)$ forces $\hat{y}=0$,
and the only feasible point of the inner maximization
problem is $y=0$. Hence,
\(
\widehat{\varphi}(x,0,\varepsilon)=x
\)
for every $x>0$. On the other hand, if $x=0$, then $\hat{y}\in K(0) =[0,1]$ and
\[
\widehat{K}(0,\hat{y},\varepsilon)
=
\{y\in[0,1]\mid y^2\leq\hat{y}^2+\varepsilon\}
=
\left[0,\min\left\{1,\sqrt{\hat{y}^2+\varepsilon}\right\}\right].
\]
Since $F(0,y)=y$, we obtain
\(
\widehat{\varphi}(0,\hat{y},\varepsilon)
=
\min\{1,\sqrt{\hat{y}^2+\varepsilon}\}.
\)
The minimum of this expression over $\hat{y}\in[0,1]$
is attained at $\hat{y}=0$ and equals $\sqrt{\varepsilon}$.
Therefore,
\[
\min_{\hat{y}\in K(x)}
\widehat{\varphi}(x,\hat{y},\varepsilon)
=
\begin{cases}
\sqrt{\varepsilon}, & x=0,\\
x, & x>0.
\end{cases}
\]

It follows that the infimum of Problem~\eqref{eq:Zeng-eps} is $0$, but this value is not attained. \hfill$\Diamond$
\end{example}

\begin{remark}\label{rem:existence}
One could avoid the nonattainment in Example~\ref{ex:nonattainment} by imposing conditions that ensure lower semicontinuity of $(x,\hat{y})\mapsto\widehat{K}(x,\hat{y},\varepsilon)$ for every $\varepsilon>0$. Under Assumptions~\ref{as:standard} and~\ref{as:convexity}, lower semicontinuity of $K$ is sufficient for this property, and hence for the existence of global minimizers of the relaxed problems; see also \cite{Svensson:2024}.
We instead retain Assumption~\ref{as:main}, which guarantees that every accumulation point of a minimizing sequence for the original pessimistic problem is a global minimizer, without requiring lower semicontinuity of $K$.
\end{remark}

%%%%%%%%%%%%%%%%%%%%%%%%%%%%%%%%%%%%%%%%
%%%%%%%%%%%%%%%%%%%%%%%%%%%%%%%%%%%%%%%%

\section{Approximation scheme}\label{sec:approximation}

The discussion in Section~\ref{sec:Preservation} suggests the following approximation scheme. Since the relaxed problems~\eqref{eq:Zeng-eps} need not have global minimizers, we work with $\delta$-minimizers:
\begin{enumerate}
    \item Choose positive sequences $\varepsilon_n\to 0$ and $\delta_n\to 0$.
    \item For each $n\in\N$, find a $\delta_n$-minimizer $(x_n,\hat{y}_n,y_n)$ of Problem~\eqref{eq:Zeng-eps} with $\varepsilon=\varepsilon_n$. Under the assumptions of Corollary~\ref{cor:MPCC-relaxation}, this is equivalent to finding a $\delta_n$-minimizer of its MPCC reformulation and retaining its primal variables.
    \item Take an accumulation point $x^*$ of $(x_n)$.
    \item Return $x^*$.
\end{enumerate}
This is a global approximation scheme, not yet a numerical algorithm. Local minimizers of the MPCC, or points with small KKT residuals, do not by themselves satisfy the $\delta_n$-minimizer requirement in Step~2. Practical implementations would also require stopping rules. We first establish that every returned point $x^*$ is a global minimizer of~\eqref{eq:Pessimistic}.

\begin{lemma}\label{lemma:widehatvarphi}
Under Assumption~\ref{as:standard}, the function $\widehat{\varphi}:\gph K\times\R_+\to\R$ has the following properties:
\begin{enumerate}
    \item $\widehat{\varphi}$ is upper semicontinuous.
    \item For every $(x,\hat{y})\in\gph K$, the map $\widehat{\varphi}(x,\hat{y},\cdot)$ is nondecreasing and continuous from the right.
\end{enumerate}
\end{lemma}
\begin{proof}
The sets $\widehat{K}(x,\hat{y},\varepsilon)$ are nonempty and compact: they contain $\hat{y}$ and are closed subsets of $Y$. Take a sequence $(x_n,\hat{y}_n,\varepsilon_n)\to(x,\hat{y},\varepsilon)$ in $\gph K\times\R_+$ and choose global maximizers $y_n\in\widehat{K}(x_n,\hat{y}_n,\varepsilon_n)$. Passing to a subsequence, we may assume that $y_n\to y\in Y$ and that $\limsup_n\widehat{\varphi}(x_n,\hat{y}_n,\varepsilon_n)=F(x,y)$. Continuity of the constraint functions gives $y\in\widehat{K}(x,\hat{y},\varepsilon)$, and hence
\(
F(x,y)\leq\widehat{\varphi}(x,\hat{y},\varepsilon).
\)
This proves upper semicontinuity.

For fixed $(x,\hat{y})$, the sets $\widehat{K}(x,\hat{y},\varepsilon)$ increase with $\varepsilon$, so $\widehat{\varphi}(x,\hat{y},\cdot)$ is nondecreasing. If $\varepsilon_n\geq\varepsilon$ and $\varepsilon_n\to\varepsilon$, then
\[
\widehat{\varphi}(x,\hat{y},\varepsilon)
\leq\liminf_n\widehat{\varphi}(x,\hat{y},\varepsilon_n)
\leq\limsup_n\widehat{\varphi}(x,\hat{y},\varepsilon_n)
\leq\widehat{\varphi}(x,\hat{y},\varepsilon).
\]
This proves continuity from the right, including at $\varepsilon=0$.
\end{proof}

For $x\in X$, $\hat{z}\in S(x)$ and $\varepsilon\geq 0$, the inclusion
\(
\widehat{K}(x,\hat{z},\varepsilon)\subset\widehat{K}(x,\hat{y},\varepsilon)
\)
 holds for every $\hat{y}\in K(x)$, and so
\begin{equation}\label{eq:min-threshold}
\begin{aligned}
\min_{\hat{y}\in K(x)}\widehat{\varphi}(x,\hat{y},\varepsilon)
&=\widehat{\varphi}(x,\hat{z},\varepsilon)\\
&=\max\left\{F(x,y)\left|\begin{array}{l}
y\in K(x),\\
f(x,y)\leq\displaystyle\min_{v\in K(x)}f(x,v)+\varepsilon
\end{array}\right.\right\}.
\end{aligned}
\end{equation}
Thus, minimizing with respect to $\hat{y}$ selects the worst value of $F$ over the $\varepsilon$-minimizers of the original follower's problem. In particular,
\[
\widehat{K}(x,\hat{z},0)=S(x),\qquad
\widehat{\varphi}(x,\hat{z},\varepsilon)\longrightarrow\varphi(x)
\quad\text{as }\varepsilon\downarrow 0.
\]

We now present the main convergence theorem. We write $\inf(\mathcal{Z}_{\varepsilon})$ for the infimum of the objective over the feasible triples of~\eqref{eq:Zeng-eps}. Under Assumption~\ref{as:standard}, this infimum is finite and the feasible set is nonempty. Consequently, $\delta$-minimizers exist for every $\delta>0$, even when the infimum is not attained.

\begin{theorem}\label{thm:main}
Consider Assumptions~\ref{as:standard} and~\ref{as:main}. Let $\varepsilon_n>0$ and $\delta_n>0$ converge to zero, and let $(x_n,\hat{y}_n,y_n)$ be feasible for Problem~\eqref{eq:Zeng-eps} with $\varepsilon=\varepsilon_n$. Suppose that these triples are $\delta_n$-minimizers, that is,
\begin{equation}\label{eq:global-approximation}
F(x_n,y_n)\leq\inf(\mathcal{Z}_{\varepsilon_n})+\delta_n.
\end{equation}
Then every accumulation point $x^*$ of $(x_n)$ is a global minimizer of Problem~\eqref{eq:Pessimistic}, and
\[
F(x_n,y_n)\longrightarrow\min_{x\in X}\varphi(x).
\]
Moreover, for any $\hat{z}^*\in S(x^*)$ and any $z^*\in\argmax_{z\in S(x^*)}F(x^*,z)$, the triple $(x^*,\hat{z}^*,z^*)$ is a global minimizer of Problem~\eqref{eq:Zeng-0}.
\end{theorem}
\begin{proof}
Assumptions~\ref{as:standard} and~\ref{as:main} ensure that Problem~\eqref{eq:Pessimistic} has a global minimizer $\bar{x}$. Fix $\hat{z}\in S(\bar{x})$. Feasibility of $(x_n,\hat{y}_n,y_n)$ gives $\varphi(x_n)\leq F(x_n,y_n)$, because every point of $S(x_n)$ is feasible for the corresponding lower-level problem. On the other hand, keeping $x=\bar{x}$ and $\hat{y}=\hat{z}$ provides a feasible comparison in~\eqref{eq:Zeng-eps}. Therefore,
\[
\min_{x\in X}\varphi(x)\leq\varphi(x_n)\leq F(x_n,y_n)
\leq\widehat{\varphi}(\bar{x},\hat{z},\varepsilon_n)+\delta_n.
\]
By Lemma~\ref{lemma:widehatvarphi}, the right-hand side converges to $\varphi(\bar{x})=\min_{x\in X}\varphi(x)$. This proves convergence of the objective values.

Let $x^*$ be any accumulation point of $(x_n)$ and pass to a subsequence, not relabeled, such that $x_n\to x^*$. Lower semicontinuity of $\varphi$ gives
\[
\varphi(x^*)\leq\liminf_n\varphi(x_n)\leq\lim_n F(x_n,y_n)
=\min_{x\in X}\varphi(x).
\]
Thus, $x^*$ solves~\eqref{eq:Pessimistic}. Finally, $\widehat{K}(x^*,\hat{z}^*,0)=S(x^*)$, so the stated choices of $\hat{z}^*$ and $z^*$ give a feasible triple of~\eqref{eq:Zeng-0} with objective value $\varphi(x^*)$. Proposition~\ref{prop:zeng-equivalence} completes the proof.
\end{proof}

% A global minimizer of~\eqref{eq:Zeng-0} can then be reconstructed as follows:
% \begin{enumerate}
%     \item Obtain a global minimizer $x^*$ of~\eqref{eq:Pessimistic} from the approximation scheme.
%     \item Choose $\hat{z}^*\in S(x^*)$ and then choose
%     \[
%     z^*\in\argmax_z\{F(x^*,z)\mid z\in K(x^*),\ f(x^*,z)\leq f(x^*,\hat{z}^*)\}.
%     \]
%     \item Return $(x^*,\hat{z}^*,z^*)$.
% \end{enumerate}
% Then $x^*$ solves~\eqref{eq:Pessimistic}, $z^*$ realizes its worst-case upper-level value, and $(x^*,\hat{z}^*,z^*)$ solves~\eqref{eq:Zeng-0}.

\section{Direct extension to robust pessimistic constraints}\label{sec:robust}

Consider the problem
\begin{equation}\label{eq:Pessimistic-Robust}
\tag{RPBP}
    \begin{array}{cl}
         \displaystyle\min_{x\in X}&\theta(x)\\
         \text{s.t.}&G(x,y)\leq 0,\quad\forall y\in S(x),
    \end{array}
\end{equation}
where $G:X\times\R^q\to\R^m$, $\theta:X\to\R$, and $S(x)$ is the solution set of~\eqref{eq:lower-level}. This is the general form of pessimistic bilevel programming with robust coupling constraints. This setting includes Problem~\eqref{eq:Pessimistic} through an epigraphical reformulation; see, e.g., \cite{Wiesemann:2013,Henke-et-al:2026:PessimisticCoupling}. We consider the following adaptation of Assumption~\ref{as:standard}.

\begin{assumption}\label{as:standard-modified} We consider the following standard assumptions:
\begin{enumerate}
    \item[(i)] $X\subset\R^p$ and $Y\subset\R^q$ are nonempty compact sets.
    \item[(ii)] The functions $G:X\times\R^q\to\R^m$ and $f,h_i,g_j:X\times\R^q\to\R$ are continuous on $X\times Y$.
    \item[(iii)] For every $x\in X$, the set $K(x)$ defined in~\eqref{eq:Feasible-lower-level} is nonempty and contained in $Y$.
    \item[(iv)] The function $\theta:X\to\R$ is lower semicontinuous.
\end{enumerate}
\end{assumption}

We can adapt the Zeng reformulation to write~\eqref{eq:Pessimistic-Robust} as an optimistic bilevel problem with coupling constraints:
\begin{equation}\label{eq:Zeng-0-Robust}
\tag{$\mathcal{Z}_{\mathrm{R}}$}
	\begin{array}{cl}
		\displaystyle\min_{x,\hat{y},(y^k)}&\theta(x)\\
		\text{s.t.}&\begin{cases}
			x\in X,\ \hat{y}\in K(x),\\
            G_k(x,y^k)\leq 0,\quad\forall k\in[m],\\
			(y^k)\text{ solves }\left\{\begin{array}{cl}
				\displaystyle\max_{(y^k)}&\displaystyle\sum_{k=1}^mG_k(x,y^k)\\
				\text{s.t.}&\begin{cases}
					y^k\in K(x),\quad&\forall k\in[m],\\
					f(x,y^k)\leq f(x,\hat{y}),\quad&\forall k\in[m].
				\end{cases}
			\end{array}\right.
		\end{cases}
	\end{array}
\end{equation}
Here $(y^k)$ denotes the tuple $(y^1,\ldots,y^m)$. For each coupling constraint $G_k$, we introduce a copy $y^k$ of the follower's variable and maximize $G_k(x,y^k)$ over $K(x)$ subject to $f(x,y^k)\leq f(x,\hat{y})$. Since these problems have independent variables, they can be combined into the single lower-level problem in~\eqref{eq:Zeng-0-Robust}.

\begin{proposition}\label{prop:robust-equivalence} Under Assumption~\ref{as:standard-modified}, Problems~\eqref{eq:Pessimistic-Robust} and~\eqref{eq:Zeng-0-Robust} are equivalent in the sense that
\[
x\text{ solves~\eqref{eq:Pessimistic-Robust}}\iff\exists\hat{y},(y^k)\text{ such that }(x,\hat{y},(y^k))\text{ solves~\eqref{eq:Zeng-0-Robust}}.
\]
\end{proposition}
\begin{proof}
A tuple $(y^k)$ solves the lower-level problem in~\eqref{eq:Zeng-0-Robust} if and only if each $y^k$ maximizes $G_k(x,\cdot)$ over $\widehat{K}(x,\hat{y})$. Indeed, its feasible set is a Cartesian product and its objective is the sum of the independent component objectives.

If $(x,\hat{y},(y^k))$ is feasible for~\eqref{eq:Zeng-0-Robust}, then $S(x)\subset\widehat{K}(x,\hat{y})$. Thus, for every $z\in S(x)$ and $k\in[m]$,
\[
G_k(x,z)\leq G_k(x,y^k)\leq 0,
\]
so $x$ is feasible for~\eqref{eq:Pessimistic-Robust}. Conversely, given a feasible $x$ for~\eqref{eq:Pessimistic-Robust}, choose $\hat{y}\in S(x)$ and $y^k\in\argmax_{y\in S(x)}G_k(x,y)$. These choices give a feasible tuple of~\eqref{eq:Zeng-0-Robust}, because $\widehat{K}(x,\hat{y})=S(x)$. The two problems therefore have the same feasible leader decisions and the same objective at each of them, proving the equivalence.
\end{proof}

The main difficulty in working with~\eqref{eq:Zeng-0-Robust} is the presence of coupling constraints. Relaxing the lower-level sublevel constraints as in~\eqref{eq:Zeng-eps} can exclude leaders that were feasible for the original problem. Instead, for $\rho>0$, we consider the penalized problem
\begin{equation}\label{eq:Zeng-0-Robust-penalized}
\tag{$\mathcal{Z}_{\mathrm{R}}^{\rho}$}
	\begin{array}{cl}
		\displaystyle\min_{x,\hat{y},(y^k)}&\displaystyle\theta(x)+\rho\sum_{k=1}^m\max\{0,G_k(x,y^k)\}\\
		\text{s.t.}&\begin{cases}
			x\in X,\ \hat{y}\in K(x),\\
			(y^k)\text{ solves }\left\{\begin{array}{cl}
				\displaystyle\max_{(y^k)}&\displaystyle\sum_{k=1}^mG_k(x,y^k)\\
				\text{s.t.}&\begin{cases}
					y^k\in K(x),\quad&\forall k\in[m],\\
					f(x,y^k)\leq f(x,\hat{y}),\quad&\forall k\in[m].
				\end{cases}
			\end{array}\right.
		\end{cases}
	\end{array}
\end{equation}

The development of the previous sections can therefore be directly applied to Problem~\eqref{eq:Zeng-0-Robust-penalized}, under the extra assumption that the functions
\begin{equation}\label{eq:InnerValueFunctionRobust}
    \begin{aligned}
        \psi_k: X&\to \R\\
        x&\mapsto \max_{y\in S(x)}G_k(x,y)
    \end{aligned}
\end{equation}
are lower semicontinuous for every $k\in [m]$. Note that this condition is also used to ensure the closedness of the feasible set of Problem~\eqref{eq:Pessimistic-Robust}, and therefore existence of solutions for the original problem.

\begin{proposition}\label{prop:penalization-convergence}
Consider Assumption~\ref{as:standard-modified}. Suppose that Problem~\eqref{eq:Pessimistic-Robust} is feasible and that the functions $\psi_k$, $k\in[m]$, defined in \eqref{eq:InnerValueFunctionRobust} are lower semicontinuous. Let $\rho_n>0$ tend to $+\infty$ and let $\delta_n\geq 0$ tend to zero. For each $n$, let $(x_n,\hat{y}_n,(y_n^k))$ be a $\delta_n$-minimizer of Problem~\eqref{eq:Zeng-0-Robust-penalized} with $\rho=\rho_n$. Then every accumulation point $x^*$ of $(x_n)$ is a global minimizer of~\eqref{eq:Pessimistic-Robust}. 
% Moreover, choosing any $\hat{z}^*\in S(x^*)$ and
% \[
% z^{k,*}\in\argmax_{z\in S(x^*)}G_k(x^*,z),\qquad k\in[m],
% \]
% gives a global minimizer $(x^*,\hat{z}^*,(z^{k,*}))$ of~\eqref{eq:Zeng-0-Robust}.
\end{proposition}
\begin{proof}
Fix any feasible point $x_0$ of~\eqref{eq:Pessimistic-Robust}. By Proposition~\ref{prop:robust-equivalence}, it can be completed to a feasible tuple of~\eqref{eq:Zeng-0-Robust}. This tuple is feasible for every penalized problem and has zero penalty. Hence,
\[
\theta(x_n)+\rho_n\sum_{k=1}^m\max\{0,G_k(x_n,y_n^k)\}
\leq\theta(x_0)+\delta_n.
\]
Since $\theta$ has a finite minimum on $X$, we obtain
\begin{align*}
0\leq\sum_{k=1}^m\max\{0,\psi_k(x_n)\}
&\leq\sum_{k=1}^m\max\{0,G_k(x_n,y_n^k)\}\\
&\leq\frac{\theta(x_0)-\min_{x\in X}\theta(x)+\delta_n}{\rho_n}
\longrightarrow 0.
\end{align*}
Let $x^*$ be an accumulation point and pass to a subsequence, not relabeled, such that $x_n\to x^*$. For every $k\in[m]$, lower semicontinuity gives
\[
\max\{0,\psi_k(x^*)\}\leq\liminf_n\max\{0,\psi_k(x_n)\}=0.
\]
Thus $x^*$ is feasible for~\eqref{eq:Pessimistic-Robust}. Moreover,
\[
\theta(x^*)\leq\liminf_n\theta(x_n)\leq\theta(x_0).
\]
Since the comparison holds for every feasible $x_0$, the point $x^*$ is a global minimizer. %Proposition~\ref{prop:robust-equivalence} and the choices in the statement give the final conclusion.
\end{proof}

The preceding arguments treat relaxation for fixed $\rho$ and penalization for $\rho\to+\infty$ separately. A simultaneous limit requires additional conditions controlling the interaction between the two parameters. Alternatively, a finite exact penalty, when available, would allow $\rho$ to be kept fixed. This analysis is beyond the scope of this note.

\paragraph{Acknowledgments:} The first and second authors where partially funded by FONDECYT grant 1251159, (ANID-Chile). The second author was partially funded by the Center of Mathematical
Modeling, CMM FB210005 BASAL (ANID-Chile).

\paragraph{Data Availability Statement:} No data was used during this research.

\paragraph{Use of IA:} GPT-6 Astra was used to revise the writing and correctness of the first draft of this work, as well as for completing the literature review. Example~\ref{ex:nonattainment} was obtained using GPT-6 Astra.

\bibliographystyle{abbrv}
\bibliography{biblioZeng}

@misc{AusselSalasSchmidt:2026,
      title={Are Single-Leader-Multi-Follower problems particular
  cases of Bilevel Programming?}, 
      author={Aussel, Didier and Salas, David and Schmidt, Martin},
      year={2026},
      howpublished={Personal Communication}
}

@article{Zeng:2020,
 author = {Zeng, Bo},
 title = {A practical scheme to compute the pessimistic bilevel optimization problem},
 fjournal = {INFORMS Journal on Computing},
 journal = {INFORMS J. Comput.},
 issn = {1091-9856},
 volume = {32},
 number = {4},
 pages = {1128--1142},
 year = {2020},
 language = {English},
 doi = {10.1287/ijoc.2019.0927},
 zbMATH = {7303827},
 Zbl = {1547.90208}
}

@article{Wiesemann:2013,
 author = {Wiesemann, Wolfram and Tsoukalas, Angelos and Kleniati, Polyxeni-Margarita and Rustem, Ber{\c{c}}},
 title = {Pessimistic bilevel optimization},
 fjournal = {SIAM Journal on Optimization},
 journal = {SIAM J. Optim.},
 issn = {1052-6234},
 volume = {23},
 number = {1},
 pages = {353--380},
 year = {2013},
 language = {English},
 doi = {10.1137/120864015}
}

@article{LamparielloEtAl:2019,
 author = {Lampariello, Lorenzo and Sagratella, Simone and Stein, Oliver},
 title = {The standard pessimistic bilevel problem},
 fjournal = {SIAM Journal on Optimization},
 journal = {SIAM J. Optim.},
 issn = {1052-6234},
 volume = {29},
 number = {2},
 pages = {1634--1656},
 year = {2019},
 language = {English},
 doi = {10.1137/18M119759X}
}

@article{LiuEtAl:2018,
  title={Pessimistic Bilevel Optimization: A Survey},
  author={June Liu and Yuxin Fan and Zhong Chen and Yue Zheng},
  year={2018},
  journal={International Journal of Computational Intelligence Systems},
  volume={11},
  issue={1},
  pages={725-736},
  issn={1875-6883},
  doi={10.2991/ijcis.11.1.56}
}

@book{solodov2010constraint,
	title={Constraint qualifications},
	author={Solodov, Mikhail V },
	year={2010},
    editors={James J. Cochran, et al.},
	PUBLISHER = {Wiley Encyclopedia of Operations Research and Management Science. James J. Cochran, et al. (editors), John Wiley \& Sons, Inc.},
    address = {New York},
}

@article {guignard1969generalized,
	AUTHOR = {Guignard, Monique},
	TITLE = {Generalized {K}uhn-{T}ucker conditions for mathematical
	programming problems in a {B}anach space},
	JOURNAL = {SIAM J. Control},
	FJOURNAL = {SIAM J. Control},
	VOLUME = {7},
	YEAR = {1969},
	PAGES = {232--241},
	MRCLASS = {90.60},
	MRNUMBER = {252042},
	MRREVIEWER = {W.\ Oettli},
}

@article{DempeDutta2012,
 author = {Dempe, S. and Dutta, J.},
 title = {Is bilevel programming a special case of a mathematical program with complementarity constraints?},
 fjournal = {Mathematical Programming. Series A. Series B},
 journal = {Math. Program.},
 issn = {0025-5610},
 volume = {131},
 number = {1-2 (A)},
 pages = {37--48},
 year = {2012},
 language = {English},
 doi = {10.1007/s10107-010-0342-1}
}

@article{LoridanMorgan:1996,
 author = {Loridan, Pierre and Morgan, Jacqueline},
 title = {Weak via strong {Stackelberg} problem: {New} results},
 fjournal = {Journal of Global Optimization},
 journal = {J. Glob. Optim.},
 issn = {0925-5001},
 volume = {8},
 number = {3},
 pages = {263--287},
 year = {1996},
 language = {English},
 doi = {10.1007/BF00121269}
}

@incollection{CarusoLignolaMorgan:2020,
 author = {Caruso, Francesco and Lignola, M. Beatrice and Morgan, Jacqueline},
 title = {Regularization and approximation methods in {Stackelberg} games and bilevel optimization},
 booktitle = {Bilevel optimization. Advances and next challenges},
 isbn = {978-3-030-52118-9; 978-3-030-52121-9; 978-3-030-52119-6},
 pages = {77--138},
 year = {2020},
 publisher = {Springer},
 address = {Cham},
 language = {English},
 doi = {10.1007/978-3-030-52119-6_4}
}

@article{Henke-et-al:2026:PessimisticCoupling,
 author = {Henke, Dorothee and Lefebvre, Henri and Schmidt, Martin and Th{\"u}rauf, Johannes},
 title = {On coupling constraints in pessimistic linear bilevel optimization},
 fjournal = {Journal of Optimization Theory and Applications},
 journal = {J. Optim. Theory Appl.},
 volume = {210},
 number = {2},
 pages = {12},
 note = {Id/No 25},
 year = {2026},
 doi = {10.1007/s10957-026-03026-x}
}

@article{Henke-et-al:2025:Coupling,
 author = {Henke, Dorothee and Lefebvre, Henri and Schmidt, Martin and Th{\"u}rauf, Johannes},
 title = {On coupling constraints in linear bilevel optimization},
 fjournal = {Optimization Letters},
 journal = {Optim. Lett.},
 volume = {19},
 number = {3},
 pages = {689--697},
 year = {2025},
 language = {English},
 doi = {10.1007/s11590-024-02156-3}
}

@article{Svensson:2024,
 author = {Svensson, Anton},
 title = {Lower semicontinuity of intersections of set-valued maps and applications on bilevel games},
 fjournal = {Set-Valued and Variational Analysis},
 journal = {Set-Valued Var. Anal.},
 volume = {32},
 number = {3},
 pages = {14},
 year = {2024},
 language = {English},
 doi = {10.1007/s11228-024-00728-0}
}

@book{Rockafellar-Wets:1998,
 author = {Rockafellar, R. Tyrrell and Wets, Roger J.-B.},
 title = {Variational analysis},
 fseries = {Grundlehren der Mathematischen Wissenschaften},
 series = {Grundlehren Math. Wiss.},
 issn = {0072-7830},
 volume = {317},
 isbn = {3-540-62772-3},
 year = {1998},
 publisher = {Springer},
 address = {Berlin},
 language = {English},
 doi = {10.1007/978-3-642-02431-3},
 zbMATH = {1110799},
 Zbl = {0888.49001}
}

@book{Dempe:2002,
 author = {Dempe, Stephan},
 title = {Foundations of bilevel programming},
 fseries = {Nonconvex Optimization and Its Applications},
 series = {Nonconvex Optim. Appl.},
 issn = {1571-568X},
 volume = {61},
 isbn = {1-4020-0631-4},
 year = {2002},
 publisher = {Kluwer Academic Publishers},
 address = {Dordrecht},
 language = {English},
 zbMATH = {1788462},
 Zbl = {1038.90097}
}

@misc{Loridan-Morgan:1989,
 author = {Loridan, P. and Morgan, J.},
 title = {{{\(\epsilon\)}}-regularized two-level optimization problems: {Approximation} and existence results},
 year = {1989},
 language = {English},
 howpublished = {Optimization, {Proc}. 5th {French}-{German} {Conf}., {Varetz}/{Fr} 1988, {Lect}. {Notes} {Math}. 1405, 99-113 (1989).},
 zbMATH = {4128662},
 Zbl = {0688.90045}
}

@article{Zemkoho1:2025,
 author = {Benchouk, Imane and Jolaoso, Lateef and Nachi, Khadra and Zemkoho, Alain},
 title = {Scholtes relaxation method for pessimistic bilevel optimization},
 fjournal = {Set-Valued and Variational Analysis},
 journal = {Set-Valued Var. Anal.},
 issn = {1877-0533},
 volume = {33},
 number = {2},
 pages = {37},
 year = {2025},
 language = {English},
 doi = {10.1007/s11228-025-00747-5},
 zbMATH = {8025030}
}

@article{Zemkoho2:2026,
 author = {Benchouk, Imane and Jolaoso, Lateef and Nachi, Khadra and Zemkoho, Alain},
 title = {Relaxation methods for pessimistic bilevel optimization},
 fjournal = {Set-Valued and Variational Analysis},
 journal = {Set-Valued Var. Anal.},
 issn = {1877-0533},
 volume = {34},
 number = {1},
 pages = {37},
 year = {2026},
 language = {English},
 doi = {10.1007/s11228-025-00787-x},
 zbMATH = {8165376},
 Zbl = {1585.90115}
}

@misc{Stein-Zemkoho:2026,
      title={A true single–level reformulation for pessimistic bilevel optimization}, 
      author={Stein, Oliver and Zemkoho, Alain},
      year={2026},
      howpublished={Optimization Online},
      url={https://optimization-online.org/?p=35995} 
}

@article{DempeEtAl:2014,
 author = {Dempe, S. and Mordukhovich, B. S. and Zemkoho, A. B.},
 title = {Necessary optimality conditions in pessimistic bilevel programming},
 fjournal = {Optimization},
 journal = {Optimization},
 issn = {0233-1934},
 volume = {63},
 number = {4},
 pages = {505--533},
 year = {2014},
 language = {English},
 doi = {10.1080/02331934.2012.696641},
 zbMATH = {6315982},
 Zbl = {1302.90206}
}

@article{GoyalZhangHe:2023,
 author = {Goyal, Akshit and Zhang, Yiling and He, Chuan},
 title = {Decision rule approaches for pessimistic bilevel linear programs under moment ambiguity with facility location applications},
 fjournal = {INFORMS Journal on Computing},
 journal = {INFORMS J. Comput.},
 volume = {35},
 number = {6},
 pages = {1342--1360},
 year = {2023},
 language = {English},
 doi = {10.1287/ijoc.2022.0168}
}

@article{CalveteEtAl:2025,
 author = {Calvete, Herminia I. and Gal{\'e}, Carmen and Hern{\'a}ndez, Aitor and Iranzo, Jos{\'e} A.},
 title = {A novel approach to pessimistic bilevel problems. An application to the rank pricing problem with ties},
 fjournal = {Optimization},
 journal = {Optimization},
 volume = {74},
 number = {12},
 pages = {2823--2856},
 year = {2025},
 language = {English},
 doi = {10.1080/02331934.2024.2388204}
}

@article{BucareyEtAl:2026,
 author = {Bucarey L{\'o}pez, V{\'i}ctor and Calder{\'o}n, Sophia and Mu{\~n}oz, Gonzalo and Semet, Fr{\'e}d{\'e}ric},
 title = {Decision-focused predictions via pessimistic bilevel optimization: complexity and algorithms},
 fjournal = {Computational Optimization and Applications},
 journal = {Comput. Optim. Appl.},
 year = {2026},
 language = {English},
 doi = {10.1007/s10589-026-00818-9}
}

@article{HoheiselKanzowSchwartz:2013,
 author = {Hoheisel, Tim and Kanzow, Christian and Schwartz, Alexandra},
 title = {Theoretical and numerical comparison of relaxation methods for mathematical programs with complementarity constraints},
 fjournal = {Mathematical Programming},
 journal = {Math. Program.},
 volume = {137},
 number = {1-2},
 pages = {257--288},
 year = {2013},
 language = {English},
 doi = {10.1007/s10107-011-0488-5}
}

@article{KanzowSchwartz:2013,
 author = {Kanzow, Christian and Schwartz, Alexandra},
 title = {A new regularization method for mathematical programs with complementarity constraints with strong convergence properties},
 fjournal = {SIAM Journal on Optimization},
 journal = {SIAM J. Optim.},
 volume = {23},
 number = {2},
 pages = {770--798},
 year = {2013},
 language = {English},
 doi = {10.1137/100802487}
}
\bigskip

\rule{5 cm}{0.5 mm} \bigskip

    Joaquín \textsc{Parra}, David \textsc{Salas}\\
    Instituto de Ciencias de la Ingeniería, Universidad de O'Higgins\\
    Rancagua, Chile\\
    E-mail: \href{mailto:joaquin.parra@pregrado.uoh.cl}{\nolinkurl{joaquin.parra@pregrado.uoh.cl}},
    \href{mailto:david.salas@uoh.cl}{\nolinkurl{david.salas@uoh.cl}}\\
\end{document}